\documentclass[11pt]{article}
\usepackage{amsfonts}
\usepackage{amsmath}
\newcommand{\oper}[2]{\newcommand{#1}{\mathop{\mathrm{#2}}\nolimits} }
\oper{\tr}{tr} \oper{\adj}{adj} \oper{\Div}{div} \oper{\ad}{ad}
\oper{\Ad}{Ad} \oper{\End}{End} \oper{\Hom}{Hom} \oper{\Aut}{Aut}
\oper{\SO}{SO} \oper{\SP}{Sp} \oper{\SU}{SU} \oper{\GL}{GL}
\oper{\T}{T} \oper{\U}{U} \oper{\id}{I} \oper{\ext}{Ext}
\oper{\rank}{rank} \oper{\diag}{Diag}

\def\bCP{{\bf CP}}

\newcommand{\CC}{\mathbb{C}}

\newcommand{\RR}{\mathbb{R}}
\newcommand{\ZZ}{\mathbb{Z}}

\newtheorem{theorem}{Theorem}
\newtheorem{corollary}[theorem]{Corollary}
\newtheorem{proposition}[theorem]{Proposition}

\newtheorem{lemma}[theorem]{Lemma}

\newcommand{\bproof}{\noindent{\it Proof: }}
\newcommand{\eproof}{\  q.~e.~d. \vspace{0.2in}}
\begin{document}

\title{Calabi-Yau Connections with Torsion\\ on Toric Bundles}
\author{D. Grantcharov \thanks{ Address:
    Department of Mathematics, San Jos\'e State University, San Jos\'e, CA 95192 , U.S.A..}
    \and G. Grantcharov \thanks{Address:
   Department of Mathematics, Florida International University, Miami, FL 33199, U.S.A..
      Partially supported by NSF DMS-0333172.}
\and Y.~S.  Poon \thanks{ Address:
    Department of Mathematics, University of California at Riverside,
    Riverside CA 92521, U.S.A..
    Partially supported by NSF DMS-0204002.}  }


\maketitle

\abstract{We find sufficient conditions for principal toric
bundles over compact K\"ahler manifolds to admit Calabi-Yau
connections with torsion as well as conditions to admit strong
K\"ahler connections with torsion. With the aid of a topological
classification, we construct such geometry on $(k-1)(S^2\times
S^4)\# k(S^3\times S^3)$ for all $k\geq 1$.

\section{Introduction}
In this article, we investigate a construction of Hermitian connections with special holonomy on Hermitian
non-K\"ahlerian manifolds. On Hermitian manifolds,
 there is a  one-parameter family of Hermitian connections
 canonically depending on the complex structure $J$ and the Riemannian metric $g$ \cite{G1}.
Among them is the Chern connection on holomorphic tangent bundles.
In this paper, we are interested in what physicists call the
K\"ahler-with-torsion connection (a.k.a. KT connection)
\cite{Strom}. It is the unique Hermitian connection whose torsion
tensor is totally skew-symmetric when 1-forms are identified to
their dual vectors with respect to the Riemannian metric. If $T$
is the torsion tensor of a KT connection, it  is characterized by
the identity \cite{G1}
\[
g(T(A, B), C)=dF(JA, JB, JC)
\]
where $F$ is the K\"ahler form; $F(A,B)=g(JA,B)$,  and $A, B, C$
are any smooth vector fields.

 As a Hermitian
connection, the holonomy of a KT connection is contained in the
unitary group $\U(n)$. If the holonomy of the KT connection is
reduced to $\SU(n)$, the Hermitian structure is said to be
Calabi-Yau with torsion (a.k.a. CYT).

Such geometry in physical context was considered first by A.
Strominger \cite{Strom} and C. Hull \cite{Hu}. More recently CYT
structures on non-K\"ahler manifolds  attracted attention as
models for string compactifications. Many examples were found
\cite{BD} \cite{BBDG} \cite{CCDLMZ} \cite{CCDL} \cite{GMW}
\cite{GP}. This led to a conjecture \cite{GIP} that {\it any}
compact complex manifold with vanishing first Chern class admits a
Hermitian metric and connection with totally skew-symmetric
torsion and (restricted) holonomy in $\SU(n)$. Counterexamples to
this conjecture appear in \cite{GF}. There are also examples of
CYT connections unstable under deformations. These two features of
CYT connections are in sharp contrast to well known moduli theory
of Calabi-Yau (K\"ahler) metric. In this paper, we shall see more
examples relevant to the moduli problem for CYT connections.

Below we begin our general construction on toric bundles on
Hermitian manifolds. Inspired by the recent results of \cite{GP}
we will focus our attention to two-dimensional bundles of torus
over compact complex surfaces. The main technical observation is
the following.

\

\noindent\bf Proposition}\hspace{0.2in}  Suppose that $X$ is a
compact K\"ahler manifold. Let the harmonic part of the Ricci form
of its K\"ahler metric be $\rho^{\mbox{har}}$. Suppose that $M$ is
a principal toric bundle with curvature $(\omega_1, \dots,
\omega_{2k})$ and that all curvature forms are harmonic type
(1,1)-forms. Then
 $M$ admits a KT connection with restricted holonomy in $\SU(n)$ if
$
\rho^{\mbox{har}}=\sum_{\ell=1}^{2k}(\Lambda\omega_\ell)\omega_\ell
$, where $\Lambda$ is a contraction with respect to the K\"ahler
form on $X$.

\

Combining the above technical observation with various algebraic
geometrical and algebraic topological results, we find a large
class of examples of compact simply-connected CYT manifolds. A
slight modification of our construction produces strong KT (a.k.a.
SKT) structures. SKT structures appear in physics literature and
refer to Hermitian structures with $dd^cF=0$  \cite{FPS}
\cite{HP}. Combining Theorem \ref{main CYT} on a construction of
CYT structures and Theorem \ref{main SKT} on a construction of SKT
structures, we establish the following observation.

\

\noindent{\bf Theorem}\hspace{.1in} For any positive integer
$k\geq 1$, the manifold $(k-1)(S^2\times S^4)\#k(S^3\times S^3)$
admits a CYT structure and a SKT structure.

\

The CYT structure and SKT structure in the above theorem do not
necessarily coincide. A CYT structure which satisfies also
$dd^cF=0$ is called strong CYT. It remains a challenge to see if
these manifolds admit strong CYT structures. On the other hand,
the existence of complex structures on these spaces is well known
\cite{OV}.

Our constructive approach to CYT structures and SKT structures is
in contrast to the obstruction theories developed by other authors
 \cite{GF} \cite{George}. It also enriches the set of examples
found in \cite{GIP}.

Although this paper focuses on constraints on KT connections, much
of its methods could be modified to construct canonical
connections subjected to similar constraints. The departure point
would be Proposition \ref{kahler}.

In the rest of this article, by a CYT connection we mean a KT
connection having {\it restricted} holonomy in $\SU(n)$. Strictly
speaking, we should have called it locally CYT. Obviously, on
simply connected manifolds such as $(k-1)(S^2\times S^4)\#
k(S^3\times S^3)$, the distinction between restricted holonomy and
holonomy disappears.

\section{Canonical connections on toric bundles}
Let $g$ be a Riemannian metric and $J$ be an integrable complex
structure such that they together form a Hermitian structure on a
manifold $M$. Let $F$ be the K\"ahler form; $F(A, B):=g(JA, B)$.
Let $d^c$ be the operator $(-1)^nJdJ$ on n-forms \cite{Besse} and
$D$ be the Levi-Civita connection of the metric $g$. Then a family
of canonical connections is given by
\begin{equation}
g(\nabla^t_AB, C)
=g(D_AB, C)+\frac{t-1}{4}(d^cF)(A,B,C)+\frac{t+1}{4}(d^cF)(A, JB, JC),
\end{equation}
where $A, B, C$ are any smooth vector fields and the real number
$t$ is a free parameter \cite{G1}. The connections $\nabla^t$ are
called \it canonical \rm connections. The connection $\nabla^{1}$
is the Chern connection on the holomorphic tangent bundle. The
connection $\nabla^{-1}$ is called the KT connection by physicists
and the Bismut connection by some mathematicians. The mathematical
features and background of these connections are articulated in
\cite{G1}. When the Hermitian metric is a K\"ahler metric, the
entire family of canonical connections collapses to a single
connection, namely the Levi-Civita connection.

In this section, we construct these connections on toric bundles
over Hermitian manifolds. We begin with a standard construction of
complex structures.

\begin{lemma}\label{toric} Suppose that $M$ is the total space of a
principal toric bundle over a Hermitian manifold $X$ with
characteristic classes of type (1,1). If the fiber is
even-dimensional, then $M$ admits an integrable complex structure
so that the projection map from $M$ to $X$ is holomorphic.
\end{lemma}
\bproof Choose a connection $(\theta_1,\theta_2,...,\theta_{2k})$
on the principal bundle $M$. Let $\pi$ denote the projection from
$M$ onto $X$. The curvature form of this connection is
$(d\theta_1,d\theta_2,...,d\theta_{2k})$. By assumption, for each
$j$ there exists (1,1)-form $\omega_j$ on $X$ such that
$d\theta_{j}=\pi^*\omega_j$.

To construct an almost complex structure $J$ on  $M$, we use the
horizontal lift of the base complex structure on the horizontal
space of the given connection. The vertical space consists of
vectors tangent to an even-dimensional torus, hence carries a
complex structure. We choose $J$ so that
$J\theta_{2j-1}=\theta_{2j}$ for $1\leq j\leq k$.

Since the fibers are complex submanifolds, if $V$ and $W$ are
vertical (1,0) vector fields, then $[V, W]$ is a vertical (1,0)
vector field. If $A^h$ and $B^h$ are horizontal lifts of (1,0)
vector fields $A$ and $B$ on $X$, $[A^h, B^h]=\omega(A, B)$. This
is equal to zero because $\omega$ is of type (1,1). Finally, if
$V$ is a vertical (1,0) vector field and $A^h$ is the horizontal
lift of a (1,0) vector field $A$ on the base manifold $X$, then
$[V, A^h]=0$ because horizontal distributions are preserved by the
action of the structure group of a principal bundle. It follows
that the complex structure on $M$ is integrable.

Since the projection from $M$ onto $X$ preserves type decomposition,
it is holomorphic.
\eproof

Suppose that $g_X$ is a Hermitian metric on the base manifold $X$
with K\"ahler form $F_X$. Consider a Hermitian metric $g_M$ on $M$
defined by
\begin{equation}\label{metric}
g_M := \pi^*g_X + \sum_{\ell=1}^{2k} (\theta_\ell\otimes
\theta_\ell).
\end{equation}
Since $J\theta_{2j-1}=\theta_{2j}$, the  K\"ahler form of the
metric $g_M$ is
\begin{equation}\label{Kahler}
F_M= \pi^*F_X+ \sum_{j=1}^k\theta_{2j-1}\wedge \theta_{2j}.
\end{equation}
Here we use the convention that $\theta_1\wedge\theta_2=\theta_1\otimes\theta_2-\theta_2\otimes\theta_1$.

Let $\Lambda$ be the contraction of differential forms on the manifold
$X$ with respect to the K\"ahler form $F_X$. If $e_1, \dots, e_{2n}$ is
a local Hermitian frame on $X$ such that $Je_{2a-1}=e_{2a}$ for $1\leq a\leq n$,
and $\omega$ is a type (1,1) form, then
\[
\Lambda\omega=\sum_{a=1}^n\omega(e_{2a-1}, e_{2a}).
\]

\begin{lemma}\label{co-differential}
If $\delta F_M$ and $\delta F_X$ are the
codifferentials of the K\"ahler forms on $M$ and $X$ respectively, then
\begin{equation}
\delta F_M = \pi^*\delta F_X  + \sum_{\ell=1}^{2k} \pi^*(\Lambda \omega_\ell)\theta_\ell
\end{equation}
\end{lemma}
\bproof Extend a local Hermitian frame $\left\{
e_{1},\ldots,e_{2n}\right\}  $ on an open subset of $X$ to a
Hermitian frame $\left\{  e_{1},\ldots
,e_{2n},t_{1},\ldots,t_{2k}\right\}  $ on an open subset of $M$
such that the vector fields $\left\{  t_{1},\ldots,t_{2k}\right\}
$ are dual to the 1-forms $\left\{
\theta_{1},\ldots,\theta_{2k}\right\}  $.

\qquad Recall that the co-differential of a tensor could be
expressed in terms
of the contraction of Levi-Civita connection $D$.%

\begin{align*}
\delta F_{M}(V)  &  =-\sum_{a=1}^{n}\left(
D_{e_{2a-1}}F_{M}\right)
(e_{2a-1},V)-\sum_{a=1}^{n}\left(  D_{e_{2a}}F_{M}\right)  (e_{2a},V)\\
&  -\sum_{j=1}^{k}\left(  D_{t_{2j-1}}F_{M}\right)
(t_{2j-1},V)-\sum _{j=1}^{k}\left(  D_{t_{2j}}F_{M}\right)
(t_{2j},V).
\end{align*}
It is a standard calculation to show that for any smooth vector
fields $A,B,$ and $C$ on a Hermitian manifold,
\[
-2(D_{A}F)(B,C)=dF(A,JB,JC)-dF(A,B,C).
\]

It follows that
\begin{align*}
\delta F_{M}(V)  &  =\sum_{a=1}^{n}dF_{M}(e_{2a-1},e_{2a},JV)+\sum_{j=1}%
^{k}dF_{M}(t_{2j-1},t_{2j},JV)\\
&  =\sum_{a=1}^{n}dF_{M}(JV,e_{2a-1},e_{2a}).
\end{align*}

Due to (\ref{Kahler}),
\begin{eqnarray}
dF_M &=& \pi^*dF_X+ \sum_{j=1}^k(d\theta_{2j-1}\wedge
\theta_{2j}-\theta_{2j-1}\wedge d\theta_{2j})
\nonumber\\
&=&\pi^*dF_X+ \sum_{j=1}^k(\pi^*\omega_{2j-1}\wedge
\theta_{2j}-\theta_{2j-1}\wedge \pi^*\omega_{2j}) \label{d-Kahler}
\end{eqnarray}
Therefore,
\begin{eqnarray} & & \delta F_M(V)=
\sum_{a=1}^{n}\pi^*dF_X(JV, e_{2a-1}, e_{2a}) \nonumber
\\&
&+\sum_{j=1}^k \sum_{a=1}^{n}\pi^*\omega_{2j-1}(e_{2a-1}, e_{2a})
\theta_{2j}(JV) -\sum_{j=1}^k
\sum_{a=1}^{n}\pi^*\omega_{2j}(e_{2a-1}, e_{2a})\theta_{2j-1}(JV)
\nonumber\\
&=&\pi^*\delta F_X (V)
-\sum_{j=1}^k(\Lambda\omega_{2j-1})(J\theta_{2j})(V)
+\sum_{j=1}^k(\Lambda\omega_{2j})(J\theta_{2j-1})(V)
\nonumber\\
&=&\pi^*\delta F_X (V)
+\sum_{j=1}^k(\Lambda\omega_{2j-1})\theta_{2j-1}(V)
+\sum_{j=1}^k(\Lambda\omega_{2j})\theta_{2j}(V).
\end{eqnarray}
So, the proof is now complete. \eproof

Hermitian manifolds with $\delta F = 0$ are called {\it balanced}
\cite{MLM}.

\begin{lemma}\label{rho1}
Let $\rho^1_X$ and $\rho^1_M$ be the Ricci forms of the Chern connections on $X$ and $M$
respectively. Then
$ \rho^1_M = \pi^* \rho^1_X.$
\end{lemma}
\bproof Let $\Theta_M$ and $\Theta_X$ be the holomorphic tangent
bundles of $X$ and $M$ respectively. Since $M$ is a holomorphic
principal toric bundle over $X$, $\Theta_M$ fits into the
following exact sequence of holomorphic vector bundles:
\[
0\to {\underline{\CC}}^k \to \Theta_M \to \pi^*\Theta_X \to 0,
\]
where ${\underline{\CC}}^k$ is the rank-k trivial holomorphic
vector bundle on $M$. It follows that the pull-back map induces a
holomorphic isomorphism between the canonical bundles $K_M$ and
$K_X$ over $M$ and $X$ respectively: $K_M=\pi^*K_X.$ On the other
hand, as a differentiable vector bundle $\Theta_M$ is isomorphic
to the direct sum ${\underline{\CC}}^k\oplus \pi^*\Theta_X$.
Therefore, the induced Hermitian metric on $K_M$ is isometric to
the pull-back metric on $K_X$. Due to the uniqueness of Chern
connection in terms of Hermitian structure and holomorphic
structure, the induced Chern connection on $K_M$ is the pull-back
of the induced Chern connection on $K_X$. Therefore, the curvature
on $K_M$ is the pull-back of the curvature of $K_X$. The same can
be said for the curvatures on the anti-canonical bundles. Since up
to a universal constant, the Ricci form is equal to the curvature
form of the Chern connection on anti-canonical bundle, the
proposition follows. \eproof

On any manifold $Y$ with a Hermitian metric $g_Y$, each canonical
connection $\nabla^t$ on the holomorphic tangent bundle induces a
connection on the anti-canonical bundle $K^{-1}_Y$. Let $R^{t}$ be
the curvature of $\nabla^t$ and $\rho_Y^t$ be the Ricci form. Then
$i\rho_Y^t$ is the curvature of the induced connection of
$\nabla^t$ on $K_Y^{-1}$. By \cite[(2.7.6)]{G1}, for any smooth
section $s$ of the anti-canonical bundle $K_Y^{-1}$ and for any
real numbers $t$ and $u$,
\begin{equation}
\nabla^t s - \nabla^u s = i \frac{t-u}{2} \delta F_Y \otimes s.
\end{equation} It follows that
\begin{equation}\label{rho}
\rho^t_Y - \rho^u_Y = \frac{t-u}{2} d\delta F_Y.
\end{equation}

\begin{proposition}\label{rhot} Let $\rho^t_M$ and $\rho^t_X$ be the Ricci forms of the canonical connections
on $M$ and $X$ respectively, then
\begin{equation}
\rho^t_M=\pi^*\rho^t_X+\frac{t-1}2\sum_{\ell=1}^{2k}d((\Lambda\omega_\ell) \theta_\ell).
\end{equation}
\end{proposition}
\bproof
Applying (\ref{rho}) on $M$ and using Lemma \ref{co-differential}, we find that
\begin{eqnarray*}
\rho^t_M-\rho^1_M
&=& \frac{t-1}2d\delta F_M
=\frac{t-1}2\left( \pi^*d\delta F_X+\sum_{\ell=1}^{2k}d\left( (\Lambda\omega_\ell)\theta_\ell
\right)\right)
\nonumber\\
&=&\pi^*(\rho_X^t-\rho_X^1)+\frac{t-1}2\sum_{\ell=1}^{2k}d\left( (\Lambda\omega_\ell)\theta_\ell
\right).
\end{eqnarray*}
Now the conclusion follows Lemma \ref{rho1}. \eproof

\begin{proposition}\label{kahler}
  Suppose that the base manifold $X$ is compact and the metric $g_X$ is K\"ahler. Let $\rho_X$ be the
Ricci form of $g_X$. If each curvature form $\omega_\ell$ is
chosen to be harmonic, then
\begin{equation}
\rho^t_M=\pi^*\left(\rho_X
+\frac{t-1}2\sum_{\ell=1}^{2k}(\Lambda\omega_\ell)\omega_\ell
\right).
\end{equation}
\end{proposition}
\bproof Every curvature form $\omega_\ell$ is closed. Up to an
addition of an exact 2-form, we may assume that $\omega_\ell$ is
harmonic. It amounts to modifying the connection 1-form
$\theta_\ell$ by the pullback of a 1-form on $X$.

Since $\omega_\ell$ is a harmonic (1,1)-form and the metric is
K\"ahler, its trace is constant \cite[2.33]{Besse}. Therefore,
\[
d\pi^*((\Lambda\omega_\ell)\theta_\ell)=\pi^*d((\Lambda\omega_\ell)\theta_\ell)
=\pi^*(\Lambda\omega_\ell)d\theta_\ell
=\pi^*((\Lambda\omega_\ell)\omega_\ell).
\]
As $g_X$ is a K\"ahler metric, all Ricci forms $\rho^t_X$ are
equal to the Ricci form $\rho_X$ of the Levi-Civita connection.
The proposition follows Lemma \ref{rhot}. \eproof

\section{CYT connections}

When the holonomy of the KT connection is contained in the special
unitary group, the KT connection is called a  CYT connection.
Locally, it is determined by the vanishing of its corresponding
Ricci form. Since the KT connection is uniquely determined by the
Hermitian structure, we address the Hermitian structure as a CYT
structure when the KT connection is CYT. In this section, we focus
on toric bundles over compact K\"ahlerian bases with various
geometrical or differential topological features.

The last proposition implies that the metric $g_M$ is CYT if the
base manifold $X$ is K\"ahler and its Ricci curvature satisfies
the following.
\begin{equation}\label{KE}
\rho_X=\sum_{\ell=1}^{2k}(\Lambda\omega_\ell)\omega_\ell.
\end{equation}
However, it is not easy to find solutions to this equation as it
requires the Ricci form of a compact K\"ahler manifold to be
harmonic. We could significantly relax the above condition by the
following observation.

\begin{lemma}
Suppose that the Ricci form of the KT connection of a Hermitian
metric $g_M$ is $dd^c$-exact on a manifold $M$ of dimension
greater than two. Then the metric $g_M$ is conformally a CYT
structure. In other words, there exists a conformal change of
$g_M$ such that the Ricci form of the induced KT connection
vanishes.
\end{lemma}
\bproof Let $\phi$ be an everywhere positive function on the
manifold $M$.  Let $\tilde{g}_M=\phi^2 g_M$ be a conformal change.
The corresponding K\"ahler forms are related by $\tilde{F}_M =
\phi^2 F_M$. The Ricci forms of the Chern connections are related
by
\begin{equation}
\tilde{\rho}_M^1 = \rho^1_M -mdd^c \log\phi
\end{equation}
where $m$ is the complex dimension of $M$ \cite[Equation
(21)]{G2}. The change of $d\delta F_M$ is
\[
d\tilde{\delta} \tilde{F}_M = d\delta F_M - 2(m-1)dd^c \log\phi.
\]
Given the universal relation among canonical connections
(\ref{rho}), we derive the relation between the KT connections of
conformally related metrics (the formula appears also in
\cite[Section 17, Lemma 1]{GIP}):
\[
{\tilde{\rho}}^{-1}_M = \rho^{-1}_M + (m-2)dd^c\log\phi.
\]
If there exists a function $\Psi$ such that
$\rho^{-1}_M=dd^c\Psi$, then
\[
{\tilde{\rho}}^{-1}_M =dd^c(\Psi + (m-2)\log\phi).
\]
Given $\Psi$, one could solve the equation $\Psi+(m-2)\log\phi=0$
for $\phi$ with $\phi(p)>0$ for every point $p$ on $M$. Therefore,
 $\tilde{\rho}_M^{-1}=0$. \eproof

\begin{proposition}\label{conformal}
Suppose that $X$ is a compact K\"ahler manifold with Ricci form
$\rho_X$. The toric bundle
 $M$ admits a CYT connection if $\rho_X^{\mbox{\rm har}}$, the harmonic
 part of $\rho_X$, satisfies the following:
\begin{equation}\label{eq conformal}
\rho_X^{\mbox{\rm
har}}=\sum_{\ell=1}^{2k}(\Lambda\omega_\ell)\omega_\ell.
\end{equation}
\end{proposition}
\bproof By
$\partial{\overline\partial}$-Lemma, there exists a function $\Phi$ such that
\[
\rho_X=\rho^{\mbox{har}}_X+dd^c\Phi.
\]
The assumption on $\rho^{\mbox{har}}_X$ and Proposition
\ref{kahler} together imply that
\[
\rho_M^{-1}=\pi^*(\rho_X-\rho^{\mbox{har}}_X)=\pi^*dd^c\Phi=dd^c\pi^*\Phi.
\]
Due to the last lemma, the metric $g_M$ is conformally equivalent to a CYT metric.
\eproof

Note that Equation (\ref{eq conformal}) has a topological
interpretaion. Since the curvature form $\omega_{\ell}$ is a
harmonic (1,1)-form, $g(\omega_{\ell},F_X) =
\Lambda(\omega_{\ell})$ is a constant. Therefore,
\[
\int_Xg_X(\omega_{\ell},F_X)dvol_X = g_X(\omega_{\ell},F_X)vol_X =
\frac{g_X(\omega_{\ell},F_X)}{n!}\int_X F_X^n.
\]
On the other hand,
\[
\int_Xg_X(\omega_{\ell},F_X)dvol_X = \int_X\omega_{\ell}\wedge\ast
F_X = \frac{1}{(n-1)!}\int_x\omega_{\ell}\wedge F_X^{n-1}.
\]
Therefore, Equation (\ref{eq conformal}) is reformulated in terms
of cohomology classes and their intersections as follows:
\[
c_1(M) = n \sum_{\ell=1}^{2k} \frac{[ \omega_{\ell} ]\cup
[F]^{n-1} }{[F]^n} [\omega_{\ell}]
\]
When the complex dimension of the manifold is equal to 2, we have
the next result.

\begin{corollary}\label{conf remark} Suppose in addition
that the base manifold $X$ is complex two-dimensional. Let $Q$ be
the intersection form of $X$. Then $M$ admits a CYT connection if
\begin{equation}\label{inter}
\rho_M^{\mbox{har}}=2\frac{Q(F_X, \omega_1)}{Q(F_X, F_X)}\omega_1
   +2\frac{Q(F_X, \omega_2)}{Q(F_X, F_X)}\omega_2.
\end{equation}
\end{corollary}

When the base manifold is K\"ahler Einstein, we could solve
Equation (\ref{KE}) directly without going through a conformal
change as in Proposition \ref{conformal}.

\begin{proposition}\label{KE base} Suppose that $X$ is compact
real $2n$-dimensional K\"ahler Einstein manifold with positive
scalar curvature. Let its scalar curvature be normalized to be
$2n^2$. Suppose that $M$ is an even-dimensional toric bundle with
curvature $(\omega_1, \dots, \omega_{2k})$ such that $\omega_1=
F_X$ and for all $2\leq \ell\leq 2k$, $\omega_\ell$ is primitive,
then $M$ admits a CYT structure.
\end{proposition}
\bproof Since $X$ is K\"ahler Einstein,
$\rho^t_X=\frac{2n^2}{2n}F_X=n F_X$ for all $t$. When
$\omega_\ell$ is primitive, $\Lambda\omega_\ell=0$. By Proposition
\ref{kahler}
\begin{equation}
\rho^{-1}_M=\pi^*(n
F_X-(\Lambda\omega_1)\omega_1)=n\pi^*(F_X-\omega_1)=0.
\end{equation}
\eproof

\begin{corollary} Let $P$ be the principal $U(1)$-bundle of the maximum root of the
anti-canonical bundle of a compact K\"ahler Einstein manifold with
positive scalar curvature, then $P\times S^1$ admits a
CYT-structure.
\end{corollary}

\begin{proposition} Let $X$ be a compact Ricci-flat K\"ahler manifold. Suppose
that $M$ is an even-dimensional toric bundle with curvature
$(\omega_1, \dots, \omega_{2k})$ such that every $\omega_\ell$ is
primitive, then there is a Hermitian metric on $M$ such that all
canonical connections are Ricci flat. In particular, it admits a
CYT structure.
\end{proposition}
\bproof Since $X$ is Ricci-flat K\"ahler,  $\rho^t_X=0$ for all
$t$. When all $\omega_\ell$ are primitive, $\Lambda\omega_\ell=0$.
By Proposition \ref{kahler} or Equation (\ref{KE}),
\begin{equation}
\rho^{t}_M=\frac{t-1}2\pi^*\left(\sum_\ell(\Lambda\omega_\ell)\omega_\ell\right)=0
\end{equation}
for all $t$.
\eproof

The condition on $\omega_\ell$ being primitive in the last
proposition is necessary as there exists an example  of real
2-dimensional holomorphic principal toric bundle over real
4-dimensional flat torus admitting no CYT connections
\cite[Theorem 4.2]{GF}.

The last proposition is applicable to K3-surfaces with Calabi-Yau
metrics. The abundance of primitive harmonic (1,1)-forms generates
a large collection of CYT structures on toric bundles on
K3-surfaces. Some explicit constructions can be found in
\cite{GP}. We shall remark on the topology of these examples and
their relation with Strominger's equations later in this article.

On the other hand,  the last two propositions could not be
extended to include K\"ahler Einstein manifolds with negative
scalar curvature as base manifolds. For instance, if $X$ is a
compact complex surface of general type, the space of holomorphic
sections of $K_X^m$ for some positive integer $m$ is at least two
dimensional. Since $K_M=\pi^*K_X$, we have a contradiction to a
vanishing theorem \cite{AI}.

\section{Examples of CYT structures}

\subsection{Product of Spheres $S^3\times S^3$.}\label{product}
The second integral cohomology of the product of two complex
projective lines $X=\bCP^1\times \bCP^1$ is generated by two
effective divisor classes $C$ and $D$ with the properties that
\[
Q(C , C)=Q(D,  D)=0, \quad Q(C,  D)=1.
\]
They are the pullback of the hyperplane class from the respective
factors onto the product space. The anti-canonical class is
$-K_X=2C+2D$. The class $C+D$ is positive as its associated map
embeds $\bCP^1\times \bCP^1$ into the complex projective 3-space
$\bCP^3$. Therefore, $F_X:=\frac12(C+D)$ is a K\"ahler class. Let
\[
\omega_1=C,
\quad
\omega_2=D.
\]
Then $ Q(F_X, F_X)=\frac12, $ and $Q(F_X,
\omega_1)=Q(F_X,\omega_2)=\frac12.$ Therefore,
\[
2\frac{Q(F_X, \omega_1)}{Q(F_X, F_X)}\omega_1
   +2\frac{Q(F_X, \omega_2)}{Q(F_X, F_X)}\omega_2
=2C+2D=-K_X.
\]
By Proposition \ref{conformal} there exists a CYT structure on the
total space of the toric bundle with curvature $(\omega_1,
\omega_2)$. Since $\omega_1$ and $\omega_2$ are the curvature of
the Hopf bundle on $\bCP^1$, the total space of the toric bundle
is simply $S^3\times S^3$. The existence of CYT structure or a
flat invariant Hermitian connection on this space is well known
\cite{FPS}.

\subsection{Toric bundles on blow-up of $\bCP^2$ twice}\label{two points}
Let $X$ be the blow-up of $\bCP^2$ at two distinct points. Let $H$
be the hyperplane class of the complex projective plane and
$E_\ell$ be the exceptional divisor of blowing-up the $\ell$-th
point on the complex projective plane. The anti-canonical class is
$3H-E_1-E_2$, and it is ample. The class $H-2E_1-E_2$ is primitive
with respect to the anti-canonical class. By taking
\begin{equation}\label{2p2}
F_X=\omega_1=3H-E_1-E_2 \quad \mbox{ and } \quad
\omega_2=H-2E_1-E_2,
\end{equation}
we solve the geometric equation in Corollary \ref{conf remark}.
Therefore, there exists a CYT structure of the total space of the
bundle whose curvature is $(\omega_1, \omega_2)$.

\subsection{Blow-up of $\bCP^2$ three to eight times}\label{KES}

Let $X$ be the blow-up of $\bCP^2$ at $k$ distinct points at
general position on the complex projective plane. Assume that
$3\leq k\leq 8$. It is well known that the anti-canonical class on
$X$ is positive and the manifold admits a K\"ahler Einstein metric
with positive scalar curvature.

Let $H$ be the hyperplane class of the complex projective plane
and $E_\ell$ be the exceptional divisor of blowing up the
$\ell$-th point. The anti-canonical class is $-K_X=3H-E_1-\cdots
-E_k.$ Let
\begin{eqnarray}
\omega_0 &=& H,
\quad
\omega_1=-K_X=3H-E_1-\cdots -E_k,
\quad
\omega_2=E_1-E_2,
\label{choice of omega}
\\
\omega_j &=& E_j,
\mbox{ for all } 3\leq j\leq k. \nonumber
\end{eqnarray}
Then  $\{\omega_0, \dots, \omega_k\}$ forms an integral basis for
$H^2(X, \ZZ)$. Let $g_X$ be a K\"ahler Einstein metric on $X$
whose K\"ahler class is equal to $\omega_1$, then $\omega_2$ is
primitive with respect to $g_X$. By Proposition \ref{KE base}
 the toric bundle with  curvature $(\omega_1, \omega_2)$
admits a CYT structure.

\subsection{Blow-ups of $\bCP^2$ many times.}\label{generic}
Next for $k\geq 9$ let $X$ be the blow-up of $\bCP^2$ at $k$
distinct points on an irreducible smooth cubic curve. Let
\begin{equation}\label{first choice}
\omega_1=4H-2\sum_{\ell=1}^4E_\ell - \sum_{\ell=5}^kE_\ell, \quad
\omega_2=-H+\sum_{\ell=1}^4E_\ell.
\end{equation}
Consider a real cohomology class $F_X = nH - \sum_{\ell=1}^k
n_\ell E_\ell$ on $X$. We now seek $n$ and $n_\ell$ for $1\leq
\ell \leq k$ such that
\begin{equation}\label{restriction}
Q(F_X,F_X)=4, \quad Q(\omega_1,F_X)=Q(\omega_2,F_X)=2,
\end{equation}
because the resulting cohomology class $F_X$ will  solve Equation
(\ref{inter}) in Corollary \ref{conf remark}. The above equations
are equivalent to the following set of equations in $n$ and
$n_\ell$.
\[
n^2-\sum_{\ell=1}^k n_\ell^2 = 4, \quad
4n-2\sum_{\ell=1}^4n_\ell-\sum_{\ell=5}^k n_\ell = 2,\quad
-n+\sum_{\ell=1}^4n_\ell = 2.
\]
To illustrate the existence of a solution, we further assume that
 $n_1=n_2=n_3=n_4$ and $n_5=\cdots =n_k$. In terms of $n$, the above system becomes
 \begin{eqnarray}
& &n_5=\cdots =n_k = \frac{2n-6}{k-4}, \quad n_1=n_2=n_3=n_4=\frac14(n+2),\\
& &(3k-28)n^2 + (112-4k)n - (20k+64)=0. \label{quadratic}
\end{eqnarray}
An elementary computation demonstrates that when $k\geq 9$ this
system of equations has a solution such that $n>3$. Note that all
$n_\ell$ are strictly positive and $n>n_\ell+n_j$ for all $\ell$
and $j$.

Next, we need to demonstrate that $F_X$ with the given $n$ and
$n_\ell$ above is also a K\"ahler class. Due to an improved
Nakai-Moishezon criteria \cite{Bu}, a class $F_X$ in
$H^{(1,1)}(X)$ is K\"ahler if the following conditions are met:
\begin{enumerate}
\item $Q(F_X, F_X)>0$. \item $Q(F_X, D)>0$ for any irreducible
curve $D$ with negative self-intersection. \item $Q(F_X, C)>0$ for
an ample divisor C.
\end{enumerate}
On $X$, the divisor class $aH-E_1-\cdots -E_k$ is ample when $a$
is sufficiently large. Therefore, the last condition is fulfilled
because $n$ is positive.

According to \cite{Fr} on the blow-up of distinct points on a
smooth cubic in $\bCP^2$, the irreducible curves with negative
self-intersections are $E_\ell$, $H-E_\ell-E_j$ with $\ell\neq j$,
and the proper transform of the cubic containing every point of
blow-up when $k\geq 10$. The latter is linearly equivalent to
$-K_X=3H-\sum_{\ell=1}^kE_\ell$. Since $\omega_1+\omega_2 = -K_X$
and $F_X$ satisfy the conditions in (\ref{restriction}), $Q(-K_X,
F_X)>0$. Therefore, $F_X$ is in the K\"ahler cone as long as the
intersection numbers of $F_X$ with $E_\ell$ and with
$H-E_j-E_\ell$ are positive. It is equivalent to the constraints $
n>n_\ell+n_j$, and $n_\ell>0$ for $1\leq \ell, j\leq k.$ Since
solutions to Equation (\ref{quadratic}) fulfill these conditions,
the corresponding $F_X$ is a K\"ahler class.

\subsection{CYT structures on $(k-1)(S^2\times S^4)\#k(S^3\times S^3)$ }

We now examine the topology of the total spaces of the toric
bundles found in the last three sections.

\begin{proposition}\label{topology}
Suppose that $X$ is a compact and simply connected manifold. Let
$M$ be the total space of a principal $T^2$-bundle over $X$ with
curvature forms $(\omega_1, \omega_2)$. Suppose that the curvature
forms are part of a set of generators $\{\omega_1, \dots,
\omega_b\}$ of $H^2(X, \ZZ)$. If there exist $\alpha, \beta$ in
$H^2(X, \ZZ)$ fulfilling the following equations on $X$,
\begin{equation}\label{intersection}
 \omega_1\wedge\alpha = \pm{\mbox{\rm vol}_X},
\quad \omega_2\wedge\alpha=0, \quad
\omega_2\wedge\beta=\pm{\mbox{\rm vol}_X}, \quad
\omega_1\wedge\beta=0,
\end{equation}
then $b=b_2(X)-2$, $H^2(M,\ZZ)=\ZZ^{b_2(X)-2}$ and the cohomology
ring of $M$ has no torsion.
\end{proposition}
\bproof There exist connection forms $(\theta_1,\theta_2)$ on $M$
with curvatures $(d\theta_1, d\theta_2)=\pi^*(\omega_1,
\omega_2)$. Let $T_x^2$ be the fiber of $M$ over a point $x$ on
the base manifold. Then  the restrictions $\theta_1|_{T^2_x}$ and
$\theta_2|_{T^2_x}$
 generate the $\ZZ$-module
$H^*(T^2_x,\ZZ)$. The standard Leray's theorem implies that the
$E_2$-terms of the Leray spectral sequence for the $T^2$-bundle
$M$ over $X$ are given in the table below \cite[15.11]{Bott}.
\[
\begin{tabular}{|c|c|c|c|c|}
\hline $\ZZ=<\theta_1\wedge\theta_2>$ & 0 & $\ZZ^b$ & 0 &
$\ZZ$ \\
\hline $\ZZ^2=<\theta_1, \theta_2>$ & 0 &
$\ZZ^2 \otimes \ZZ^b$ & 0 & $\ZZ^2$ \\
\hline $\ZZ$ & 0 & $\ZZ^b=<\omega_1,\omega_2,...,\omega_b>$
& 0 & $\ZZ= <\mbox{vol}_X>$ \\
\hline
\end{tabular}.
\]
Next we calculate the $E_3$-terms. The map $d_2: E_2^{0,1}
\rightarrow E_2^{2,0}$ is given by $ \theta_1,\theta_2, \mapsto
\omega_1,\omega_2.$ It is an injection and therefore
$E_3^{0,1}=0$. Since $d_2(E_2^{2,0})=0$, $E_3^{2,0}=\ZZ^{b-2}\cong
<\omega_3,\omega_4,...,\omega_b>$.

 Similarly, the map
$d_2:E^{0,2}_2 \rightarrow E^{2,1}_2$ is given by $
\theta_1\wedge\theta_2 \mapsto
\omega_1\wedge\theta_2-\theta_1\wedge\omega_2$ so it is injective.
Therefore, $E_3^{0,2}=0$. The map $d_2:E_2^{2,1}\rightarrow
E^{4,0}_2$ is given by $ \theta_i\wedge\omega_j\mapsto
\omega_i\wedge\omega_j. $ It is surjective because
$d_2(\theta_1\wedge\alpha)=\pm\mbox{vol}_X$. Therefore,
$E_3^{2,1}= \ZZ^{2b-2}$.

Finally, the map $d_2:E^{2,2}_2 \rightarrow E_2^{4,1}$ is given by
\[
\theta_1\wedge\theta_2\wedge\omega \mapsto
\omega_1\wedge\theta_2\wedge\omega-\theta_1\wedge\omega_2\wedge\omega
\]
for any $\omega\in H^2(X,\mathbb{Z})$. In particular,
\[
d_2(\theta_1\wedge\theta_2\wedge\alpha)=\pm
\mbox{vol}_X\wedge\theta_2 \quad \mbox{and} \quad
d_2(\theta_1\wedge\theta_2\wedge\beta)=\pm
\mbox{vol}_X\wedge\theta_1.
\]
Therefore, the restriction $d_2$ on the $E^{2,2}_2$-term  is
surjective. It follows that $E_3^{2,2}=\ZZ^{b-2}$. With the other
$E_3$ terms easily computed, the above computation yields the
table of $E_3$-terms below.
\[
\begin{tabular}{|c|c|c|c|c|}
\hline
0 & 0 & $\ZZ^{b-2}$ & 0 & $\ZZ$ \\
\hline
0 & 0& $ \ZZ^{2b-2}$ & 0 & 0 \\
\hline
$\ZZ$ & 0 & $\ZZ^{b-2}$ & 0 & 0 \\
\hline
\end{tabular}
\]
It follows that the spectral sequence degenerates at the
$E_3$-level and
\[
H^2(M, \ZZ)\cong E_3^{2,0}\oplus E_3^{1,1}\oplus E_3^{0,2}
\cong E_3^{2,0}\cong <\omega_3,\omega_4,...,\omega_b> .
\]
\eproof

\begin{theorem}\label{main CYT}
For every positive integer $k\geq 1$, the manifold
$(k-1)(S^2\times S^4)\#k(S^3\times S^3)$ admits a CYT structure.

Alternatively, any 6-dimensional compact simply-connected spin manifold
with torsion free cohomology and free $S^1$-action admits a CYT structure.
\end{theorem}
\bproof We have seen a CYT structure on $S^3\times S^3$ in a
previous section.

For $k\geq 2$, let  $X_k$ be the blow-up  of $\bCP^2$ at $k$
distinct points on an irreducible smooth cubic. Let $M_k$ be the
total space of the toric bundles over $X_k$ obtained in Sections
\ref{two points}, \ref{KES}, or \ref{generic}. Since $M_k$ admits
a CYT structure, $c_1(M_k)=0$. In particular, the second
Stiefel-Whitney class vanishes and $M_k$ is a spin manifold.

We now examine the cohomology of the space $M_k$ through the last
proposition. When $k=2$, let $\omega_1$ and $\omega_2$ be given as
in (\ref{2p2}) and let
\[
\alpha=H+E_1-3E_2, \quad \beta=E_1-E_2.
\]
When $3\leq k\leq 8$, let $\omega_1$ and $\omega_2$ be given as in
(\ref{choice of omega}) and let
\[
\alpha=E_k, \quad \beta=E_1-E_k.
\]
When $k\geq 9$, let $\omega_1$ and $\omega_2$ be  given as in
(\ref{first choice}) and let
\[
\alpha=E_k, \quad \beta=H-E_5-E_6-E_7-E_8.
\]
 Then the set of data $\{\omega_1,
\omega_2, \alpha, \beta\}$ on respective manifolds satisfies the
hypothesis of Proposition \ref{topology}. In particular,
$b_2(M_k)=k-1$, for each $k\geq 2$.

Since $M_k$ is a toric bundle, it admits a free $S^1$-action. By
\cite{GL}, a compact smooth simply-connected spin 6-manifold with
torsion-free cohomology, $b_2(M_k) = k-1$ and free $S^1$-action is
diffeomorphic to $(k-1)(S^2\times S^4)\# k(S^3\times S^3)$.
Therefore, we complete the proof of this theorem if we demonstrate
that the space $M_k$ is simply connected.

When $3\leq k\leq 8$, $Q(\omega_1, \alpha)=1$ and $Q(\omega_2,
\alpha)=0$, and the restriction of the bundle $M_k$ onto the
2-sphere representing the homotopy class of $\alpha$ is the
projection from $S^3\times S^1$ onto $S^2$ via the Hopf fibration
$S^3\to S^2$. In particular,  the map $\pi_2(M_k)\to \pi_1(T^2)$
in the homotopy sequence of the fibration from  $M_k$ onto $X_k$
sends $\alpha$ to a generator of $\pi_1(T^2)$ \cite{Steenrod}.
Similarly, the Poincar\'e dual of $\beta$ is represented
topologically by a smooth 2-sphere. As $Q(\omega_1, \beta)=0$ and
$Q(\omega_2, \beta)=-1$, the map $\pi_2(M_k)\to \pi_1(T^2)$ sends
the Poincar\'e dual of $\beta$ to a different generator of
$\pi_1(T^2)$. Therefore, the map $\pi_2(M_k)\to \pi_1(T^2)$ is
surjective. Since $X_k$ is simply connected, by the homotopy
sequence of the fibration $M_k \to X_k$, $M_k$ is simply
connected.

When $k\geq 9$, a similar analysis shows that the map $\pi_2(M_k)
\to \pi_1(T^2)$ sends the Poincar\'e dual of $\alpha$ and $\beta$
onto the generators of $\pi_1(T^2)$. Hence, $M_k$ is simply
connected.

When $k=2$, the Poincar\'e dual of $\beta$ is an embedded
2-sphere. Since $Q(\omega_1, \beta)=0$ and $Q(\omega_2, \beta)=1$,
the restriction of the bundle $M_2\to X_2$ onto the Poincar\'e
dual of $\alpha$ is the Hopf fibration $S^1\times S^3 \to S^2$.
The map $\pi_2(M_2)$ $\to$ $\pi_1(T^2)$ sends $\beta$ to a
generator in $\pi_1(T^2)$. The Poincar\'e dual of
$\gamma=H-E_1-E_2$ is also an embedded sphere. Since $Q(\omega_1,
\gamma)=1$, under the map $\pi_2(M_2)\to \pi_1(T^2)$ the images of
$\gamma$ and $\beta$ form a set of generators for $\pi_1(T^2)$.
\eproof

\section{SKT connections}
A KT connection is strong (a.k.a. SKT) if its torsion is a \it
closed \rm three-form. It is equivalent to require the K\"ahler
form to be $dd^c$-closed. Such structures recently appeared in the
theory of generalized K\"ahler geometry \cite{Gualtieri}
\cite{Hitchin}. In real six dimension,  a Hermitian metric with
strong KT connection is an astheno-K\"ahler metric \cite{JY}.
Results on SKT structures on nilmanifolds could be found in
\cite{FPS}.

To construct SKT connections, we return to the general set-up
leading to Lemma \ref{co-differential}. Since the projection map
$\pi$ is holomorphic and the curvature forms $(\omega_{2j-1},
\omega_{2j})$ are type (1,1), given Equation (\ref{d-Kahler}) we
have
\begin{eqnarray*}
d^cF_M &=& JdF_M
   = J\pi^*dF_X+\sum_{j=1}^k(\pi^*\omega_{2j-1}\wedge J\theta_{2j}-J\theta_{2j-1}\wedge \pi^*\omega_{2j})
\nonumber\\
 &=& \pi^*d^cF_X-\sum_{j=1}^k(\pi^*\omega_{2j-1}\wedge \theta_{2j-1}+\theta_{2j}\wedge \pi^*\omega_{2j}).
\end{eqnarray*}
Therefore,
$
dd^cF_M = \pi^*dd^cF_X-\sum_{j=1}^k \pi^*(\omega_{2j-1}\wedge
\omega_{2j-1}+\omega_{2j}\wedge \omega_{2j}).
$

\begin{proposition} Suppose that a Hermitian structure on a toric bundle $M$ over a Hermitian manifold $X$
 is given as in {\rm Equation (\ref{metric})}. Its
KT connection is strong if and only if
\begin{equation}\label{strong}
\sum_{j=1}^k (\omega_{2j-1}\wedge \omega_{2j-1}+\omega_{2j}\wedge
\omega_{2j})=dd^cF_X.
\end{equation}
In particular, suppose that $X$ is a compact complex K\"ahler
surface, $Q$ is the intersection form on $X$, and  $M$ is a real
two-dimensional toric bundle over $X$. If the KT connection on $M$
is strong, then
\begin{equation}\label{strong2}
Q(\omega_{1}, \omega_{1})+Q(\omega_{2}, \omega_{2})=0.
\end{equation}
\end{proposition}

Note that unlike the trace, the square of a harmonic form is not
harmonic so (\ref{strong2}) is not equivalent to (\ref{strong}).
By Hodge-Riemann bilinear relations, the intersection form on
primitive type (1,1) classes on compact complex surfaces is
negative definite \cite{Griffiths}. There is little chance of
using our construction here to produce strong CYT structures on
2-toric bundles over K3-surfaces. However, on $S^2\times S^2$,
when $\omega_1=C$ and $\omega_2=D$ as given in Section
\ref{product}, they solve the equation (\ref{strong}). Therefore,
the CYT-structure on $S^3\times S^3$ is strong, which is a well
known fact.

It is known that on most rational surfaces a product of harmonic
forms is not harmonic and every harmonic anti-self-dual 2-form has
at least one zero \cite{Kot}.  So we must use a non-K\"ahler
metric on $X$. With this observation in mind we are ready to prove
the following:

\begin{theorem}\label{main SKT}
For every positive integer $k\geq 1$, the manifold
$(k-1)(S^2\times S^4)\# k(S^3\times S^3)$ admits a strong KT
structure.

Alternatively, any 6-dimensional compact simply-connected spin
manifold with torsion free cohomology and free $S^1$-action admits
a SKT structure.
\end{theorem}
\bproof
 Consider a blow-up of $\bCP^2$ at $k$ ($k\geq 2$) points on a smooth irreducible
 cubic.
Choose two arbitrary
 closed forms $\omega_1$ and $\omega_2$ satisfying
 (\ref{strong2}). By the $dd^c$-lemma $\omega_1\wedge\omega_1 +
 \omega_2\wedge\omega_2 = dd^c\alpha$ for some real (1,1)-form
 $\alpha$.
  We can choose a $dd^c$-closed (1,1)-form $\beta$ (e.g.
 any K\"ahler form multiplied by an appropriate constant), such
 that $$min_{p\in X}(min_{||Y|| = 1}\beta_p(Y,JY)) > -min_{p\in X}(min_{||Y|| =
 1}\alpha_p(Y,JY))$$
Then the form $\alpha + \beta$ is positive definite everywhere and
defines a Hermitian metric on X, which will produce SKT metric on
$M$.
 \eproof

\section{Remarks}
\subsection{Non-uniqueness}\label{non-unique}
 The construction on bundles over the blow-ups of $\bCP^2$ in Section \ref{generic}
 could be used to
produce apparently different CYT structures on $(k-1)(S^2\times
S^4)\# k(S^3\times S^3)$ using  different K\"ahler classes on the
same base manifold. For example when $k\geq 11$, we may choose
\[
\omega_1=4H-2(E_1 + E_2) - \sum_{\ell=3}^kE_\ell, \quad
\omega_2=-H+E_1+E_2,
\] and then solve
 the equations $ Q(\omega_1,F_X)=Q(\omega_2,F_X)=2$,
$Q(F_X,F_X)=4$ for $n, a, b$ in $F_X=nH-a(E_1+E_2)-b(E_3+\cdots
+E_k)$.

It is also possible to use topologically different base manifolds
and toric bundles to produce CYT structures on the same real
six-dimensional manifold. For instance, let the base manifold $X$
be a Kummer surface. It admits sixteen  smooth rational curves
$C_i$ with $Q(C_i,C_j)=-2\delta_{ij}$. Due to Piateckii-Shapiro
and Shafarevich's description of the cohomology ring of $X$
\cite[pp.568-571]{P-S},
 if we choose $\omega_1=C_1\pm C_2$ and $\omega_2=C_3\pm C_4$ (with arbitrary signs),
 there are elements $\alpha$ and $\beta$ in $H^2(X,\ZZ)$ satisfying all conditions
in Proposition  \ref{topology}. This proposition enables us to
identify the total space $M$ of the  toric bundle with $(\omega_1,
\omega_2)$ as curvature forms of $20(S^2\times S^4)\#21(S^3\times
S^3)$. Since the canonical bundle of $M$ is the pullback of the
one on $M$, it is holomorphically trivial.

Moreover the K\"ahler cone of $X$ in $H^2(X,\RR)$ is one of the
chambers with walls $\Pi_i = \{ E\in H^2(X,\ZZ):Q(C_i,E)=0\}$.
Then we can choose a Ricci-flat K\"ahler metric $F_X$ on $X$ such
that $Q(C_i,F_X)=\pm 1$ for $i=1,2,3,4$. For this choice we have
$Q(F_X,\omega_j)=0$ for $j=1,2$ after fixing the signs in the
definition of $\omega_i$. As a consequence we obtain a balanced
CYT structure on $20(S^2\times S^4)\#21(S^3\times S^3)$ which
implicitly appears in \cite{GP}. It is also  half-flat in the
terminology of \cite{Simon}.  The structure is not strong, which
is in accordance to the result in \cite{George} that a strong and
balanced CYT structure on compact manifold is K\"ahler.

\subsection{Relation with Strominger's equations}
Our construction on the connected sums of $S^2\times S^4$ with
$S^3\times S^3$ is related to a set of Strominger's equations in
string theory. In \cite{Strom} Strominger analyzes heterotic
superstring background with spacetime supersymmetry. His model can
be translated to our situation in the following terms: First we
need a conformally balanced CYT manifold with holomorphic
(3,0)-form of constant norm. The manifold is endowed with an
auxiliary semistable bundle with Hermitian-Einstein connection $A$
with curvature $F_A$. The last and most restrictive equation in
\cite{Strom} is \begin{equation}\label{last} dH =
\alpha'({\mbox{Tr}} R\wedge R - {\mbox{tr}}F_A\wedge F_A).
\end{equation}
Here "${\mbox{Tr}}$" and "${\mbox{tr}}$" are the traces in the
tangent bundle and the auxiliary bundle respectively, $R$ is the
curvature of any metric connection $\nabla$ and $\alpha'$ is a
positive constant. Solutions to the Strominger's equations with
the choice of $R$ being the curvature of the Chern connection
$\nabla^1$ have recently been found by Fu and Yau \cite{yau}.
 The Hermitian metric has
K\"ahler form as in (2) with $F_X$ being conformally K\"ahler.
With this data, they solve the system proposed by Strominger for
an unknown conformal factor on a K3-surface as base space. Since
the anomalies can be cancelled for $any$ choice of metric
connection, it is important progress towards a realistic string
theory \cite{hep-th/0604137}. However the requirement that the
connection $\nabla$ preserves both worldsheet conformal invariance
and spacetime supersymmetry leads the connection $\nabla$ to be
equal to $D- \frac{1}{2} H$, where H is the torsion of the KT
(Bismut) connection. i.e. $\nabla^{-1}=D+\frac12 H$ \cite{Hu}.
Therefore, the term $R$ in Equation (\ref{last}) is the curvature
of the connection  $D- \frac{1}{2} H$. It is an open question
whether such a connection exists on compact manifolds.

\subsection{Orbifolds}
In this article, we focus on toric bundles over smooth complex
manifolds. However, most of the local geometric considerations
could be extended to toric bundles as  V-bundles over  orbifolds.
For example, suppose
 that the base space $X$ is a K\"ahler Einstein orbifold with positive scalar curvature
\cite{DK}. Let $P$ be the principal $U(1)$-bundle of the maximal
root of the anti-canonical bundle. Our construction shows that
$S^1\times P$ carries a CYT structure. We refer the readers to
\cite{BG} for an analysis of the geometric and topological
consideration of $P$.

\subsection{Other geometric structures}
The spaces $M=(k-1)(S^2\times S^4)\# k(S^3\times S^3)$ are
$S^1$-bundles over  the Sasakian 5-manifolds  $ (k-1)(S^2\times
S^3)$. Since $(k-1)(S^2\times S^3)$ admits a Ricci-positive metric
for any $k\geq 2$ \cite{SY} \cite{BG2}, a result of
B\'erard-Bergery \cite{BB} implies that
  the manifolds $M$ admit
 Riemannian metrics with positive Ricci curvature.
 It would be interesting to see to what extent the CYT metric,
the SKT metric and the positive Ricci curvature metrics on $M$ are
related.

\

\noindent{\bf Acknowledgement.}\  We are very grateful to Prof. S.
T. Yau for his interest in our work. We are also in debt to the
very knowledgable and thoughtful referee(s) for the extremely
useful suggestions. Finally we want to thank J.Edward for the help
in editing.

\end{document}